\documentclass[12pt]{article}
\usepackage{latexsym}
\usepackage{theorem}
\usepackage{graphicx}
\usepackage{amsmath}
\usepackage{color}
\usepackage{xcolor}
\usepackage{amsfonts}
\usepackage{natbib}
\usepackage{soul}
\usepackage{hyperref}

\theorembodyfont{\rmfamily}
\newtheorem{theorem}{Theorem}
\newtheorem{conjecture}[theorem]{Conjecture}
\newtheorem{lemma}[theorem]{Lemma}
\newtheorem{proposition}[theorem]{Proposition}

\theoremstyle{break}

\title{A Search Game on a Hypergraph with Booby Traps}

\author{Thomas Lidbetter\thanks{Department of Management Science and Information Systems, Rutgers Business School, Newark, NJ 07102, tlidbetter@business.rutgers.edu} \and Kyle Y. Lin\thanks{Operations Research Department, Naval Postgraduate School, Monterey, CA 93943, kylin@nps.edu}}

\providecommand{\keywords}[1]{\textbf{\textbf{Keywords:}} #1}

\begin{document}

\maketitle

\begin{abstract}
\noindent
A set of $n$ boxes, located on the vertices of a hypergraph $G$, contain known but different rewards. A {\em Searcher} opens all the boxes in some hyperedge of $G$ with the objective of collecting the maximum possible total reward.
Some of the boxes, however, are booby trapped.
If the Searcher opens a booby trapped box, the search ends and she loses all her collected rewards. We assume the number $k$ of booby traps is known, and we model the problem as a zero-sum game between the maximizing Searcher and a minimizing {\em Hider}, where the Hider chooses $k$ boxes to booby trap and the Searcher opens all the boxes in some hyperedge. The payoff is the total reward collected by the Searcher.
This model could reflect a military operation in which a drone gathers intelligence from guarded locations, and a booby trapped box being opened corresponds to the drone being destroyed or incapacitated.
It could also model a machine scheduling problem, in which rewards are obtained from successfully processing jobs but the machine may crash.
We solve the game when $G$ is a $1$-uniform hypergraph (the hyperedges are all singletons), so the Searcher can open just 1 box.
When $G$ is the complete hypergraph (containing all possible hyperedges), we solve the game in a few cases: (1) same reward in each box, (2) $k=1$, and (3) $n=4$ and $k=2$.
The solutions to these few cases indicate that a general simple, closed form solution to the game appears unlikely.
\end{abstract}


\newpage

\section{Introduction}
\label{sec:introduction}

We consider the following game between a Hider and a Searcher. There is a set $[n] \equiv \{1,\ldots,n\}$ of boxes, with box $i$ containing a reward of $r_i \ge 0$, for $i \in [n]$.
We also make a standing assumption that, without loss of generality, $r_1 \geq \cdots \geq r_n$. The boxes are identified with the vertices of a hypergraph $G$. The Hider sets booby traps in $k$ of the boxes, where $1 \le k \le n-1$, so his strategy set is ${[n]^{(k)} \equiv \{H \subset [n]: |H|=k \}}$. The Searcher chooses a subset $S \subset [n]$ of boxes to search, where $S$ is the hyperedge of a hypergraph $G$ with vertices $V$ and hyperedges $E \subset 2^V$.

If the Hider plays $H$ and the Searcher plays $S$, the payoff $R(S,H)$ is given by 
\[   
R(S,H) = 
     \begin{cases}
       r(S), & \text{if } H \cap S = \emptyset, \\
       0, & \text{otherwise,}
     \end{cases}
\]
where $r(S) \equiv \sum_{i \in S} r_i$ is the sum of the rewards in $S$. In other words, the Searcher keeps the sum of all the rewards in the boxes she opens unless one or more of them is booby trapped, in which case, she gets nothing.
If the Searcher uses a mixed strategy $p$ (that is, a probability distribution over subsets $S \subset [n]$) and the Hider uses a mixed strategy $q$ (a probability distribution over subsets $H \in [n]^{(k)}$), we write the expected payoff as $R(p,q)$.
We also write $R(p, H)$ and $R(S, q)$ if one player uses a pure strategy while the other player uses a mixed strategy.

This game could be an appropriate model for a military scenario in which a drone is used to gather intelligence at several locations, and $r_i$ is the expected value of the intelligence gathered at location $i$. A known number $k$ of the locations are guarded, and flying the drone near these locations would result in its incapacitation.
Alternatively, the Searcher may be collecting rewards in the form of stolen weapons or drugs from locations at which capture is possible, or the Searcher could be a burglar stealing valuable possessions from houses in a neighborhood, some of which are monitored by security cameras.
The graph structure could correspond to geographical constraints. 
The case of the complete hypergraph, where $E=2^V$, corresponds to no constraints on the Searcher's choice of subset. The case where $E$ is $1$-uniform, so that every hyperedge consists of a single vertex, corresponds to the Searcher being limited to searching only one location. If $E$ is $2$-uniform, so that $G$ is a graph, the Searcher must choose locations corresponding to the endpoints of an edge of the graph.

The game could also model a scheduling problem in which there are $n$ jobs with utilities $r_i$ which are obtained from a successful execution of job $i$.
For example, jobs may correspond to computer programs.
A total of $k$ of the programs are bugged, and each bug will crash the machine so that all data is lost. The objective is to find a subset of jobs to run that maximizes the worst-case expected utility, assuming Nature chooses which $k$ jobs are bugged.

This work lies in the field of {\em search games}, as discussed in \cite{AlpernGal}, \cite{Gal2011}, and \cite{Hohzaki}. 
Search games involving objects hidden in boxes have previously been considered in~\cite{Lidbetter} and~\cite{LidbetterLin}.
In these works, the objective of the Searcher is to {\em minimize} a total cost of finding a given number of hidden objects. \cite{Agnetis} consider a machine scheduling problem in which rewards are collected from processing jobs and the machine may crash, similarly to our problem. 
But in their setting, each job will independently cause the machine to fail with a given probability.

Since this is a zero-sum game, it could be solved by standard linear programming methods, but this approach would be inefficient for large $k$, or if the hypergraph has a large number of hyperedges. In this work, we concentrate on two special cases of the game, with the aim of finding concise, closed-form solutions.
We first solve the case where $G$ is a 1-uniform hypergraph in Section~\ref{sec:null}.
In Section~\ref{sec:complete}, we consider the complete hypergraph, and solve the game for three special cases: (1) same reward in each box, (2) $k=1$, and (3) $n=4$, $k=2$.
We also give some general bounds, and make a conjecture on the form of the optimal solution.
Finally, we offer concluding remarks in Section~\ref{sec:conclusion}.

\section{The game on a 1-uniform hypergraph}
\label{sec:null}
We begin with the special case that $G$ is a $1$-uniform hypergraph, so that every hyperedge is a singleton (though every singleton may not be a hyperedge). In other words, the Searcher can open only 1 box, and her strategy set is simply some subset $A$ of the set $[n]$ of vertices.
If the Searcher is restricted to boxes in $A$, then any Hider strategy that does not hide all $k$ booby traps in $A$ is (weakly) dominated by another Hider strategy that does.
Hence, without loss of generality, we may assume that $G$ is the complete $1$-uniform hypergraph whose hyperedges are {\em all} the singletons.
A mixed strategy for the Searcher is a probability vector $\mathbf x \in \mathbb R^n$ with $\sum_{j=1}^n x_j =1, x_j \ge 0$ for all $j$.

We first obtain a class of lower bounds on the value of the game, by defining a Searcher strategy for every subset of boxes. 
\begin{lemma} \label{lem:null}
For a subset $A \subseteq [n]$ of boxes with $|A| \geq k$, let the Searcher strategy $\mathbf x \equiv \mathbf x^A$ be given by
\[
x^A_j = 
\begin{cases}
\lambda(A)/r_j, & \text{ if } j \in A, \\
0, & \text{ otherwise,}
\end{cases}
\]
where $\lambda(A) = \left( \sum_{i \in A} 1/r_i \right)^{-1}$.
The strategy $\mathbf{x}^A$ guarantees an expected payoff of at least $(|A|-k)\lambda(A)$.
\end{lemma}
\textit{Proof.}
The expected payoff of the Searcher strategy $\mathbf{x}^A$ against the Hider's strategy $H$ is
\[
R(\mathbf{x}^A, H) = \sum_{j \in A-H} x^A_j \, r_j = |A- H| \, \lambda(A) \geq (|A|-  |H|) \, \lambda(A) = (|A|-  k) \, \lambda(A),
\]
where the lower bound is obtained when $H \subseteq A$.
\hfill $\Box$

Recall that $r_1 \geq r_2 \geq \cdots \geq r_n$.
If the Searcher is restricted to choosing a strategy of the form described in Lemma~\ref{lem:null}, for some $|A|=t \ge k$, then it is clear that the subset maximizing $(|A|-  k) \, \lambda(A)$ is $[t] = \{1, 2, \ldots, t\}$.
For $t=k, k+1, \ldots, n$, define
\begin{equation}
V(t) \equiv (t-k) \lambda ([t]),
\label{eq:V(t)}
\end{equation}
which is the expected payoff guaranteed by choosing $A = [t]$ in Lemma~\ref{lem:null}.
Our main result is that, when $G$ is the complete 1-uniform hypergraph, the value of the game is $\max_{t=k,\ldots,n} V(t)$.
For example, if $n=3$ and $k=1$ with $(r_1,r_2,r_3)=(10,10,1)$, then $V(1)=0$, $V(2) = 5$, and $V(3) = 5/3$, so the value of the game is $V(2)=5$ and the Searcher opens either box 1 or box 2 each with probability 0.5.
Intuitively, if the rewards in different boxes are lopsided, then it is better for the Searcher to avoid those boxes with the lowest rewards altogether.
We need a lemma before presenting the theorem.

\begin{lemma}
For $t \geq k+1$, the two inequalities $V(t) \geq V(t-1)$ and $r_t \geq  V(t)$ are equivalent, where $V(t)$ is defined in (\ref{eq:V(t)}).
\label{le:r_t}
\end{lemma}
\textit{Proof.}
The first inequality is equivalent to
\[
\frac{t-k}{\frac{1}{r_1} + \cdots + \frac{1}{r_t}} \geq \frac{t-1-k}{\frac{1}{r_1} + \cdots + \frac{1}{r_{t-1}}} 
\]
Multiplying both denominators on both sides and cancelling common terms yields
\[
\frac{1}{r_1} + \cdots + \frac{1}{r_t} \geq (t-k) \frac{1}{r_t},
\]
which is equivalent to $r_t \geq V(t)$, thus completing the proof.
\hfill $\Box$

\bigskip

\begin{theorem}
\label{th:null}
Consider the search game with $n$ boxes and $k$ booby traps played on the complete $1$-uniform hypergraph.
Define
\[
t^* \equiv \arg \max_{t = k,\ldots,n} V(t),
\]
where $V(t)$ is defined in \eqref{eq:V(t)}.
The strategy $\mathbf x^{[t^*]}$ described in Lemma~\ref{lem:null} is optimal for the Searcher. 
For the Hider, any strategy that distributes the $k$ booby traps among the boxes in $[t^*]$ in such a way that box $j \in [t^*]$ contains a booby trap with probability
\[
y_j \equiv 1-  \frac{V(t^*)}{r_j}
\]
is optimal.
The value of the game is $V(t^*)$.
\end{theorem}
\textit{Proof.}
By Lemma~\ref{lem:null}, the Searcher guarantees an expected payoff at least $V(t^*)$ by using the strategy $\mathbf x^{[t^*]}$, so $V(t^*)$ is a lower bound for the value of the game. 

To show that $V(t^*)$ is also an upper bound for the value of the game, first note that $t^* \ge k+1$, since $V(k) = 0 < V(k+1)$.
In addition, by definition of $t^*$, we have that $V(t^*-1) \le V(t^*)$, which is equivalent to $V(t^*) \le r_{t^*}$ by Lemma~\ref{le:r_t}, so $y_j \in [0,1]$ for $j \in [t^*]$.
One can also verify that $\sum_{j =1}^{t^*} y_j = k$.

If the Hider's strategy has the property described in the theorem, then the expected payoff against any Searcher strategy $j \in [t^*]$ is $r_j (1-y_j) = V(t^*)$ and
the expected payoff against any Searcher strategy $j \notin [t^*]$ is $r_j$.
By definition of $t^*$, we have that $V(t^*) \ge V(t^*+1)$, which is equivalent to $V(t^* +1) \ge r_{t^*+1}$ from the proof in Lemma~\ref{le:r_t}.
Combining two inequalities yields that $V(t^*) \geq V(t^*+1) \geq r_{t^*+1}$.
In other words, opening box $j \notin [t^*]$ results in payoff $r_j \leq r_{t^*+1} \leq V(t^*)$.
Consequently, $V(t^*)$ is an upper bound for the value of the game, which completes the proof.
\hfill $\Box$

\bigskip

There are many Hider strategies that will give rise to the property required in Theorem~\ref{th:null}; that is, the Hider distributes $k$ booby traps in $t^*$ boxes in such a way that box $j \in [t^*]$ contains a booby trap with probability $y_j \in [0,1]$, where $\sum_{j=1}^{t^*} y_j = k$.
One way to implement such a Hider strategy can be found in Definition 2.1 in \cite{Gal-Casas}.
Partition the interval $[0, k]$ into subintervals of lengths $y_1, \ldots, y_{t^*}$.
Generate $\theta$ from the uniform distribution in $[0,1]$ and select the $k$ boxes corresponding to the $k$ subintervals containing the points $\theta, \theta+1, \ldots, \theta+(k-1)$.
By construction, the Hider will choose exactly $k$ boxes to put booby traps, and box $i$ will contain a booby trap with probability $y_i$, for $i \in [t^*]$.

%


In the special case where all the rewards are equal, we have  $V(t) = (t-k)/t$, which is maximized at $V(n) = (n-k)/n$. 
The Searcher's optimal strategy is to open each box with probability $1/n$, and any Hider strategy that puts a booby trap in each box with the same probability $k/n$ is optimal---such as choosing every subset of $k$ boxes with probability $1/{n \choose k}$.

In another special case when $k=n-1$, we have  $V(n-1)=0$, so the value of the game is $V(n) = (n- (n-1)) \lambda([n]) = ( \sum_{j =1}^n 1/r_j )^{-1}$.
The Searcher's optimal strategy is to open box $j$ with probability $\lambda([n])/ r_j$, for $j \in [n]$.
The Hider's optimal strategy needs to put a booby trap in box $j$ with probability $1- (n - (n-1)) \lambda([n])/ r_j$.
Because the Hider has $n-1$ booby traps, the only strategy that meets this requirement is for the Hider to leave box $j$ free of booby trap with probability  $\lambda([n])/ r_j$, for $j \in [n]$.

\section{The game on the complete hypergraph}
\label{sec:complete}
This section concerns the extreme case where $G$ is the complete (non-uniform) hypergraph, so that a Searcher strategy is any $S \subset [n]$.
Note that if $k=n-1$, the Searcher should open only 1 box, so the structure of the hypergraph becomes irrelevant; the solution presented in Section~\ref{sec:null} is also optimal.

For the case of complete hypergraph, we present the solution to the three special cases: (1) equal rewards in each box; (2) $k=1$, and (3) $n=4$, $k=2$.
We then give some general bounds on the value of the game, and make a conjecture on the optimal solution based on our findings.



\subsection{The case with equal rewards} \label{sec:equal}
We begin our analysis with the special case where all the rewards are equal, which we set to $1$ without loss of generality.

\begin{theorem} \label{thm:eq}
Consider the search game on the complete hypergraph with $r_j=1$ for $j \in [n]$, so this game is characterized by only the number of boxes $n$ and the number of booby traps $k$.
The Hider's optimal strategy is to choose some $H \in [n]^{(k)}$ uniformly at random.
The Searcher's optimal strategy is to open $m^* = \lceil  \frac{n-k}{k+1} \rceil$ boxes at random.
In particular, $m^*=1$ if $k \geq \frac{n-1}{2}$.
The value of the game is given by
\[
U(n, k)  \equiv \frac{{n-m^* \choose k} m^*}{{n \choose k}}.
\]
\end{theorem}
\textit{Proof.}
By symmetry, it is optimal for the Hider to choose uniformly at random between all his pure strategies. 

Because each box contains the same reward, the Searcher's decision reduces to the number of boxes she opens.
Write $F(m)$ for the expected reward when the Searcher opens $m$ boxes at random, and the booby trap is located in some arbitrary set of $k$ boxes.
We calculate $F(m)$ by considering the Searcher's $m$ boxes to be fixed and supposing that a randomly chosen set of $k$ boxes are booby trapped. The Searcher gets a reward of $m$ if none of the boxes she has chosen are booby trapped; otherwise she gets nothing. Hence,
\begin{align}
F(m) = \frac{{n-m \choose k} m}{{n \choose k}}. \label{eq:F}
\end{align}
The ratio $F(m+1)/F(m)$ is given by
\[
\frac{F(m+1)}{F(m)} = \frac{(n-k-m)(m+1)}{(n-m)m}.
\]
Therefore, $F(m+1) \le F(m)$ if and only if
\[
m \ge \frac{n-k}{k+1}.
\]
It follows that $F(m)$ is maximized at $m^* = \lceil (n-k)/(k+1) \rceil$.
The value of the game is $F(m^*)$, as given in the statement of the theorem.
\hfill $\Box$

\medskip

Note that if $\theta \equiv k/n \le 1/2$ is held constant, and $n$ and $k$ tend to $\infty$, then the optimal search strategy in the limit is to open $m^*=\lceil (1-\theta)/\theta \rceil$ boxes, which is independent of $n$.
The same result is obtained independently in Example 2.1c in \cite{ross} with a dynamic programming formulation.
Writing $F(m)$ as
\[
F(m) = \left( \frac{n-k}{n} \right) \left( \frac{n-k-1}{n-1} \right) \cdots \left( \frac{n-k-m+1}{n-m+1} \right) m,
\]
we can verify that the value of the game in the limit is
\[
\lim_{n \rightarrow \infty} F(m^*) =m^* (1-\theta)^{m^*}.
\]
One can interpret $(1-\theta)^{m^*}$ as the probability that none of the $m^*$ boxes opened by the Hider contains a booby trap in the limit as $n \rightarrow \infty$.

Suppose now that $k$ is held constant and let $n \rightarrow \infty$.
In the limit, the optimal number of boxes to open tends to infinity, and so does the value of the game.
To calculate the proportion of the total reward $n$ the Searcher can obtain, we write out the probability ${n-m^* \choose k} / {n \choose k}$ from~(\ref{eq:F}) that none of the Searcher's boxes are booby trapped as
\begin{equation}
\prod_{i=0}^{k-1}  \left( 1- \frac{m^*}{n-i} \right),
\label{eq:success_prob}
\end{equation}
where $m^* = \lceil (n-k)/(k+1) \rceil$.
Because 
\[
\frac{n-k}{k+1} \leq m^* < \frac{n-k}{k+1} + 1 = \frac{n+1}{k+1},
\]
the probability in \eqref{eq:success_prob} satisfies the bounds 
\[
\prod_{i=0}^{k-1} \left( 1- \frac{\frac{n+1}{n-i}}{k+1} \right) < \prod_{i=0}^{k-1}  \left( 1- \frac{m^*}{n-i} \right) \leq \prod_{i=0}^{k-1}  \left( 1- \frac{\frac{n-k}{n-i}}{k+1} \right).
\]
Since the upper bound and the lower bound approach to the same limit as $n \rightarrow \infty$, we can conclude that
\[
\lim_{n \rightarrow \infty} \prod_{i=0}^{k-1}  \left( 1- \frac{m^*}{n-i} \right) = \left(1 - \frac{1}{k+1} \right)^k.
\]
Hence, in the limit as $n \rightarrow \infty$, the ratio of the value of the game to the total reward $n$ is
\begin{align}
\lim_{n \rightarrow \infty}  \frac{U(n,k)}{n} = \frac{1}{k+1} \left(1- \frac{1}{k+1} \right)^k. \label{eq1}
\end{align}

\subsection{The case with $k=1$ booby trap}
We now consider the special case in which the Hider has only $k=1$ booby trap.
Recall that a Searcher's pure strategy is $S \subset [n]$.
In order to present an optimal strategy for the Searcher, define $S^* \subset [n]$ to be a subset of boxes that minimizes $|r(S) - r(\bar{S})|$, where $\bar{S}$ denotes the complement of $S$. 

We state and prove optimal strategies for the game in the case $k=1$. Let $R_0 = \sum_{i=1}^n r_i$.

\begin{theorem} \label{theorem:k=1}
Consider the search game on the complete hypergraph with $k=1$. Let $S^* \subset [n]$ be a subset of boxes that minimizes $|r(S) - r(\bar{S})|$. It is optimal for the Searcher to choose $S^*$ with probability 
\[
p(S^*) = \frac{r(\bar{S}^*)}{R_0};
\]
otherwise choose $\bar{S}^*$. It is optimal for the Hider to put the booby trap in box $i$ with probability $q_i = r_i/R_0$, for $i \in [n]$. The value $V$ of the game is
\[
V = \frac{r(S^*) \, r(\bar{S}^*)}{R_0}.
\]
\end{theorem}
\textit{Proof.}
Suppose the Searcher uses the strategy $p$ and that the booby trap is in some box $j$. 
If $j \in S^*$, the expected payoff is
\[
p(\bar{S}^*) \, r(\bar{S}^*) = \frac{r(S^*) \, r(\bar{S}^*)}{R_0}.
\]
Similarly, if $j \in \bar{S}^*$, the expected payoff is the same.
Therefore, $V \ge r(S^*) r(\bar{S}^*)/R_0$.

On the other hand, suppose the Hider uses the strategy $q$.
If the Searcher opens some subset $S$ of boxes, then the expected payoff is
\[
r(S) \sum_{i \in \bar{S}} q_i  = \frac{r(\bar{S}) \, r(S)}{R_0}.
\]
The numerator in the preceding is equal to
\[
r(S) \, r(\bar{S}) = r(S)(R_0 - r(S)) = - \left(r(S) - \frac{R_0}{2} \right)^2 + \frac{R_0^2}{4} = -\frac{(r(S) - r(\bar{S}))^2}{4} + \frac{R_0^2}{4},
\]
which is maximized by taking $S=S^*$ by definition of $S^*$.
In other words, the Hider's strategy $q$ guarantees that the expected payoff (for the Searcher) is at most $r(S^*) r(\bar{S}^*)/R_0$, so $V \le  r(S^*) r(\bar{S}^*)/R_0$.
The result follows.
\hfill $\Box$

\bigskip

In the case that the rewards are integers, the problem of finding such a subset $S^*$ to minimize $|r(S) - r(\bar{S})|$ is the optimization version of the {\em number partitioning problem}, which is the problem of deciding whether a multiset of positive integers can be partitioned into two sets such that the sum of the integers in each set is equal.
This problem is NP-hard, so that finding the value of the search game with $k=1$ is also NP-hard.
There are, however, efficient algorithms to solve the problem in practice \citep{Korf}.


Note that the value of the game for $k=1$ is $R_0/4$, if and only if the boxes can be partitioned into two subsets of equal total reward. It is tempting to conjecture that in general, the value of the game is $R_0/(k+1)^2$, if and only if the boxes can be partitioned into $k+1$ subsets of equal total rewards.
This conjecture, however, is not true, as can be seen from the simple example with $n=6$ and $k=2$ when all the rewards are equal to $1$.
By Theorem~\ref{thm:eq}, the value of the game is $4/5$, but $R_0/(k+1)^2 = 6/3^2$. 
Nevertheless, the quantity $R_0/(k+1)^2$ is a lower bound for the value of the game, because the Searcher can choose one of the $k+1$ subsets uniformly at random, and receive an expected payoff of $R_0/(k+1)$ with probability at least $1/(k+1)$.

\subsection{The case with $n=4$ boxes and $k=2$ booby traps}
This section presents the solution to the game with $n=4$ boxes and $k=2$ booby traps. The Hider chooses two boxes to place the booby traps, so he has ${4 \choose 2} = 6$ pure strategies.
The Searcher would want to open at most $n-k = 4-2=2$ boxes, so she has 10 viable pure strategies, including ${4 \choose 1} =4$ pure strategies that open just 1 box, and ${4 \choose 2} = 6$ pure strategies that open 2 boxes.
While one can compute the value $V$ and optimal strategy of each player by a linear program, we will show that the optimal mixed strategy for the Searcher is one of the following three types:
\begin{enumerate}
	\item
	Strategy A involves 4 active pure strategies: $\{1\}$, $\{2\}$, $\{3\}$, $\{4\}$.
	Specifically, the Searcher opens just 1 box, and chooses box $i$ with probability 
	\[
	p_i = \frac{1/r_i}{1/r_1+1/r_2+1/r_3+1/r_4}, \qquad i=1,2,3,4.
	\]
	Regardless of which two boxes contain booby traps, strategy A produces the same expected payoff
	\begin{equation}
	V_A \equiv \frac{2}{\frac{1}{r_1} + \frac{1}{r_2} + \frac{1}{r_3} + \frac{1}{r_4}}.
	\label{eq:V_A}
	\end{equation}
	Intuitively, strategy A works well if $r_1, r_2, r_3, r_4$ are comparable.
	
	\item
	Strategy B involves 3 active pure strategies: $\{1\}$, $\{2\}$, and $\{3,4\}$.
	Specifically, the Searcher opens box $i$, for $i=1,2$, with probability
	\[
	p_i = \frac{1/r_i}{1/r_1 + 1/r_2 + 1/(r_3+r_4)}, \qquad i=1,2,
	\]
	or opens both boxes 3 and 4 with probability 
	\[
	p_{34} = \frac{1/(r_3+r_4)}{1/r_1 + 1/r_2 + 1/(r_3+r_4)}.
	\]
	Regardless of which two boxes contain booby traps, strategy B guarantees an expected payoff at least
	\begin{equation}
	V_B \equiv \frac{1}{\frac{1}{r_1} + \frac{1}{r_2} + \frac{1}{r_3 + r_4}}.
	\label{eq:V_B}
	\end{equation}
	Intuitively, strategy B works well if $r_3 + r_4$ is comparable to  $r_1$ and $r_2$.
	
	\item
	Strategy C involves 6 pure strategies: $\{1\}$, $\{2\}$, $\{3\}$, $\{1,4\}$, $\{2,4\}$, $\{3,4\}$.
	Specifically, the Searcher opens box $i$, for $i=1,2,3$, with probability
	\[
	p_i = \frac{1/r_i}{\frac{1}{r_1} + \frac{1}{r_2} + \frac{1}{r_3} + \frac{1}{r_1+r_4} + \frac{1}{r_2+r_4} + \frac{1}{r_3+r_4}}.
	\]
	or opens both boxes $i$ and 4, for $i=1,2,3$, with probability
	\[
	p_{i4} = \frac{1/(r_i+r_4)}{\frac{1}{r_1} + \frac{1}{r_2} + \frac{1}{r_3} + \frac{1}{r_1+r_4} + \frac{1}{r_2+r_4} + \frac{1}{r_3+r_4}}.
	\]
	Regardless of which two boxes contain booty traps, strategy C produces the same expected payoff
	\begin{equation}
	V_C \equiv \frac{2}{\frac{1}{r_1} + \frac{1}{r_2} + \frac{1}{r_3} + \frac{1}{r_1+r_4} + \frac{1}{r_2+r_4} + \frac{1}{r_3+r_4}}.
	\label{eq:V_C}
	\end{equation}
	Intuitively, strategy C works well if $r_4$ is much smaller than the reward in each of the other three boxes.
\end{enumerate}

The main result in this section is the following theorem.

\begin{theorem}
\label{th:n4k2}
One of the three strategies A, B, C is optimal for the Searcher.
The value of the game is
\[
V = \max \{ V_A, V_B, V_C\},
\]
where $V_A$, $V_B$, and $V_C$ are defined in \eqref{eq:V_A}, \eqref{eq:V_B}, and \eqref{eq:V_C}, respectively.
\end{theorem}

The proof of this theorem is lengthy, and we will present the three cases separately.
Before doing so, we first offer some discussion to shed light on these three strategies.
With some algebra, one can see that $V_A \geq V_B$ if and only if
\begin{equation}
\frac{1}{r_1} + \frac{1}{r_2} \geq \frac{1}{r_3} + \frac{1}{r_4} - \frac{2}{r_3+r_4};
\label{eq:AB}
\end{equation}
and $V_A \geq V_C$ if and only if 
\begin{equation}
\frac{1}{r_4} \leq \frac{1}{r_1+r_4} + \frac{1}{r_2+r_4} + \frac{1}{r_3+r_4};
\label{eq:AC}
\end{equation}
and $V_B \geq V_C$ if and only if
\begin{equation}
\frac{1}{r_1} + \frac{1}{r_2} + \frac{1}{r_3+r_4} \leq \frac{1}{r_3} + \frac{1}{r_1+r_4} + \frac{1}{r_2+r_4}.
\label{eq:BC}
\end{equation}

Strategy A treats each box equally.
For Strategy A to work well, $r_1$ and $r_2$ cannot be too large compared with $r_3$ and $r_4$ (seen in \eqref{eq:AB}), and $r_4$ cannot be too small (seen in \eqref{eq:AC}).
In other words, the four rewards need to be somewhat comparable.
Strategy B combines the two boxes with smaller rewards together, and treats the problem as if there were only 3 boxes.
For Strategy B to work well, $r_3$ and $r_4$ need to be substantially smaller than $r_1$ and $r_2$ (seen in \eqref{eq:AB}), and $r_3$ needs to be somewhat closer to $r_4$ rather than to $r_2$ (seen in \eqref{eq:BC}).
Strategy C treats box 4---the one with the smallest reward---as a small add-on to one of the other three boxes.
For Strategy C to work well, $r_4$ needs to be small enough (seen in \eqref{eq:AC}), and $r_1, r_2, r_3$ need to be somewhat close together (seen in \eqref{eq:BC}).

We next present the proof of Theorem~\ref{th:n4k2} in three sections, starting with the easiest case.
The challenge in each of the three proofs is to show that the Hider has a mixed strategy to guarantee the payoff to be no more than the corresponding payoff guaranteed by the Searcher's mixed strategy.

\subsubsection{Optimality of Strategy C}

\begin{theorem}
Strategy C is optimal for the Searcher and the value of the game is $V_C$ if and only if
\begin{equation}
\frac{1}{r_1+r_4} + \frac{1}{r_2+r_4} + \frac{1}{r_3+r_4} \leq \frac{1}{r_4},
\label{eq:C>=A}
\end{equation}
and
\begin{equation}
\frac{1}{r_1} + \frac{1}{r_2} + \frac{1}{r_3+r_4} \geq \frac{1}{r_3} + \frac{1}{r_1+r_4} + \frac{1}{r_2+r_4}.
\label{eq:C>=B}
\end{equation}
\end{theorem}
\textit{Proof.}
If strategy C is optimal for the Searcher, then $V_C \geq V_A$, which is equivalent to \eqref{eq:C>=A}, and $V_C \geq V_B$, which is equivalent to \eqref{eq:C>=B}.
Therefore, \eqref{eq:C>=A} and \eqref{eq:C>=B} are necessary conditions.

We next prove \eqref{eq:C>=A} and \eqref{eq:C>=B} are sufficient conditions.
Since the Searcher can use strategy C to guarantee an expected payoff $V_C$, it remains to show that the Hider has a mixed strategy to guarantee an expected payoff no more than $V_C$.
Let $q_{ij}$ denote the probability that the Hider hides the 2 booby traps in boxes $i$ and $j$, and let
\begin{align*}
q_{12} &= \frac{V_C}{r_3+r_4},  \qquad q_{13} = \frac{V_C}{r_2+r_4},  \qquad q_{23} = \frac{V_C}{r_1+r_4},
\\
q_{14} &= 1 - \frac{V_C}{r_2 + r_4} - \frac{V_C}{r_3 + r_4} - \frac{V_C}{r_1}, \\
q_{24} &= 1 - \frac{V_C}{r_1 + r_4} - \frac{V_C}{r_3 + r_4} - \frac{V_C}{r_2}, \\
q_{34} &= 1 - \frac{V_C}{r_1 + r_4} - \frac{V_C}{r_2 + r_4} - \frac{V_C}{r_3}.
\end{align*}

First, we show that the preceding is indeed a legitimate mixed strategy for the Hider.
Using the definition in \eqref{eq:V_C}, one can verity that $\sum_{1 \leq i < j \leq 4} q_{ij} =1$.
In addition, $0 \leq q_{23} \leq q_{13} \leq q_{12} \leq 1$ and $q_{34} \leq q_{24} \leq q_{14} \leq 1$, because $r_1 \geq r_2 \geq r_3 \geq r_4$.
Finally, we see that $q_{34} \geq 0$, due to \eqref{eq:C>=B}.

Next, we show that the Hider guarantees an expected payoff no more than $V_C$ regardless of what the Searcher does.
Consider 4 cases.
\begin{enumerate}
\item
If the Searcher opens $\{1,4\}$, then the expected payoff is
\[
(r_1+r_4) q_{23} = V_C.
\]
A similar argument leads to the same conclusion if the Searcher opens $\{2,4\}$ or $\{3,4\}$.
\item
If the Searcher opens \{1\}, then the expected payoff is
\[
r_1 (q_{23} + q_{24} + q_{34} ) = V_C.
\]
A similar argument leads to the same conclusion if the Searcher opens \{2\} or \{3\}.

\item
If the Searcher opens \{4\}, then the expected payoff is
\[
r_4 (q_{12}+q_{23}+q_{13})  = r_4 \left( \frac{1}{r_1+r_4} + \frac{1}{r_2+r_4} + \frac{1}{r_3+r_4} \right) V_C \leq  V_C,
\]
where the inequality follows from \eqref{eq:C>=A}.
\item
If the Searcher opens $\{1,2\}$, then the expected payoff is
\[
(r_1+r_2) q_{34} = (r_1+r_2) \left( 
\frac{1}{r_1} + \frac{1}{r_2} + \frac{1}{r_3+r_4} - \frac{1}{r_3} - \frac{1}{r_1+r_4} - \frac{1}{r_2+r_4} \right) \frac{V_C}{2}
\]
To show that the preceding is no more than $V_C$, compute
\begin{align*}
& (r_1+r_2) \left( 
 \frac{1}{r_1} + \frac{1}{r_2} + \frac{1}{r_3+r_4} - \frac{1}{r_3} - \frac{1}{r_1+r_4} - \frac{1}{r_2+r_4} \right) \\
& = (r_1+r_2) \left ( \frac{r_4}{r_1(r_1+r_4)} + \frac{r_4}{r_2(r_2+r_4)} - \frac{r_4}{r_3(r_3+r_4)} \right) \\
& = \frac{r_4}{r_1+r_4} + \frac{r_4}{r_2+r_4} + \left( \frac{r_1 r_4}{r_2(r_2+r_4)} - \frac{r_1 r_4}{r_3(r_3+r_4)} \right) + \left(\frac{r_2 r_4}{r_1(r_1+r_4)} - \frac{r_2 r_4}{r_3(r_3+r_4)}\right) \\
& \leq \frac{r_4}{r_1+r_4} + \frac{r_4}{r_2+r_4} + \frac{r_1 r_4}{r_1(r_1+r_4)} + \frac{r_2 r_4}{r_2(r_2+r_4)}  \\
&= 2  \left(\frac{r_4}{r_1+r_4} + \frac{r_4}{r_2+r_4} \right) \\
& \leq 2 \left( 1 - \frac{r_4}{r_3+r_4} \right) \leq 2,
\end{align*}
where the first inequality follows from $q_{14} \geq 0$ and $q_{24} \geq 0$, and the second inequality follows from \eqref{eq:C>=A}.
A similar argument leads to the same conclusion if the Searcher opens $\{2,3\}$ or $\{1,3\}$.
\end{enumerate}
The proof is complete.
\hfill $\Box$


\subsubsection{Optimality of Strategy B}

\begin{theorem}
Strategy B is optimal for the Searcher and the value of the game is $V_B$ if and only if
\begin{equation}
\frac{1}{r_1} + \frac{1}{r_2} \leq \frac{1}{r_3} + \frac{1}{r_4} - \frac{2}{r_3+r_4}.
\label{eq:B>=A}
\end{equation}
and
\begin{equation}
\frac{1}{r_1} + \frac{1}{r_2} + \frac{1}{r_3+r_4} \leq \frac{1}{r_3} + \frac{1}{r_1+r_4} + \frac{1}{r_2+r_4}.
\label{eq:B>=C}
\end{equation}
\end{theorem}
\textit{Proof.}
If strategy B is optimal for the Searcher, then $V_B \geq V_A$, which is equivalent to \eqref{eq:B>=A}, and $V_B \geq V_C$, which is equivalent to \eqref{eq:B>=C}.
Therefore, \eqref{eq:B>=A} and \eqref{eq:B>=C} are necessary conditions.

We next prove \eqref{eq:B>=A} and \eqref{eq:B>=C} are sufficient conditions.
Since the Searcher can use strategy B to guarantee an expected payoff at least $V_B$, it remains to show that the Hider has a mixed strategy to guarantee an expected payoff no more than $V_B$.
Let $q_{ij}$ denote the probability that the Hider hides the 2 booby traps in boxes $i$ and $j$, and require
\begin{align}
q_{12} &= \frac{V_B}{r_3+r_4}, \label{eq:12} \\
q_{34} &= 0, \label{eq:34} \\
q_{23} + q_{24} + q_{34} &= \frac{V_B}{r_1}, \label{eq:23+24} \\
q_{13} + q_{14} + q_{34} &= \frac{V_B}{r_2}. \label{eq:13+14}
\end{align}

These constraints ensure that $\sum_{1 \leq i < j \leq 4} q_{ij} = 1$, and guarantee an expected payoff no more than $V_B$ if the Searcher uses pure strategies $\{3, 4\}$, $\{1, 2\}$, \{1\}, and \{2\}.

The Hider also needs to ensure an expected payoff no more than $V_B$ if the Searcher uses either $\{3\}$ or $\{4\}$, so we need to require
\begin{align}
q_{12} + q_{14} + q_{24} &\leq \frac{V_B}{r_3}, \label{eq:12+14+24} \\
q_{12} + q_{13} + q_{23} &\leq \frac{V_B}{r_4}; \label{eq:12+13+23}
\end{align}
and if the Searcher uses $\{1,3\}$, $\{2,3\}$, $\{2,4\}$, or $\{1,4\}$, so we also need to require
\begin{align}
q_{24} &\leq \frac{V_B}{r_1+r_3}. \label{eq:24} \\
q_{14} &\leq \frac{V_B}{r_2+r_3}, \label{eq:14} \\
q_{13} &\leq \frac{V_B}{r_2+r_4}, \label{eq:13} \\
q_{23} &\leq \frac{V_B}{r_1+r_4}. \label{eq:23}
\end{align}
To complete the proof, we need to show that there exists a feasible nonnegative solution to $q_{ij}$, $1 \leq i < j \leq 4$ subject to the constraints in \eqref{eq:12} through \eqref{eq:23}.





To proceed, write
\begin{equation}
\frac{q_{13}}{V_B} = x, \qquad  \frac{q_{23}}{V_B} = y,
\label{eq:13,23}
\end{equation}
and use \eqref{eq:34} in \eqref{eq:23+24} and \eqref{eq:13+14} to obtain
\begin{equation}
\frac{q_{24}}{V_B} = \frac{1}{r_1} - y, \qquad
\frac{q_{14}}{V_B} = \frac{1}{r_2} - x.
\label{eq:14,24}
\end{equation}
To ensure $q_{13}, q_{23}, q_{14}, q_{24} \geq 0$, we need that
\begin{equation}
0 \leq x \leq \frac{1}{r_2}, \qquad 0 \leq y \leq \frac{1}{r_1}.
\label{eq:x,y}
\end{equation}

Next, subsitute \eqref{eq:13,23} and \eqref{eq:14,24} into \eqref{eq:12+14+24}--\eqref{eq:23} to rewrite the 6 inequalities constraints in terms of $x$ and $y$.
Constraints \eqref{eq:12+14+24} and \eqref{eq:12+13+23} together become
\begin{equation}
\frac{1}{r_1} + \frac{1}{r_2} - \frac{1}{r_3} + \frac{1}{r_3+r_4} \leq x + y \leq \frac{1}{r_4} - \frac{1}{r_3+r_4}.
\label{eq:x+y}
\end{equation}
Constraints \eqref{eq:14} and \eqref{eq:13} together become
\begin{equation}
\frac{1}{r_2} - \frac{1}{r_2 + r_3} \leq x \leq \frac{1}{r_2+r_4},
\label{eq:x}
\end{equation}
and constraints \eqref{eq:24} and \eqref{eq:23} together become
\begin{equation}
\frac{1}{r_1} - \frac{1}{r_1 + r_3} \leq y \leq \frac{1}{r_1+r_4}.
\label{eq:y}
\end{equation}
Because constraints \eqref{eq:x} and \eqref{eq:y} make constraint \eqref{eq:x,y} redundant, it remains to show that there exists a feasible solution to $x$ and $y$ subject to constraints \eqref{eq:x+y}, \eqref{eq:x}, and \eqref{eq:y}.

First, note that in each of \eqref{eq:x+y}, \eqref{eq:x}, and \eqref{eq:y}, the unknown's upper bound is greater than or equal to its lower bound.
The feasibility of $x+y$ in \eqref{eq:x+y} follows directly from \eqref{eq:B>=A}.
The feasibility of $x$ in \eqref{eq:x} follows from $r_2 \geq r_3 \geq r_4$, and the feasibility of \eqref{eq:y} follows from $r_1 \geq r_3 \geq r_4$.

To complete the proof, we need to show that the sum between the upper bound (lower bound, respectively) of $x$ in \eqref{eq:x} and the upper bound (lower bound, respectively) of $y$ in \eqref{eq:y} is greater than or equal to the lower bound (upper bound, respectively) of $x+y$ in \eqref{eq:x+y}.

The first claim follows directly from \eqref{eq:B>=C}.
The second claim states that
\[
\frac{1}{r_2} - \frac{1}{r_2 + r_3} + \frac{1}{r_1} - \frac{1}{r_1 + r_3} \leq  \frac{1}{r_4} - \frac{1}{r_3+r_4}.
\]
To prove it, start with the left-hand side to obtain
\begin{align*}
r_3 \left( \frac{1}{r_2(r_2 + r_3)} + \frac{1}{r_1(r_1 + r_3)} \right) &\leq  r_3 \left( \frac{1}{r_2(r_2 + r_4)} + \frac{1}{r_1(r_1 + r_4)} \right) \\ & \leq r_3 \left(\frac{1}{r_3(r_3+r_4)} \right) \\ & \leq r_3 \left(\frac{1}{r_4(r_3+r_4)} \right) \\
&= \frac{1}{r_4} - \frac{1}{r_3+r_4},
\end{align*}
where the first and third inequalities are due to $r_3 \geq r_4$, and the second inequality is due to \eqref{eq:B>=C}.
Consequently, we have proved that there exists a feasible solution to $x$ and $y$ that satisfy constraints \eqref{eq:x+y}, \eqref{eq:x}, and \eqref{eq:y}.
In other words, we have proved that there exists a feasible solution to $q_{13}, q_{14}, q_{23}, q_{24}$ that satisfy the constraints in \eqref{eq:23+24} through \eqref{eq:23}.
Therefore, we have shown that the Hider has a mixed strategy that guarantees the Searcher no more $V_B$, which completes the proof.
\hfill $\Box$

\subsubsection{Optimality of Strategy A}
We begin with two lemmas.

\begin{lemma}
\label{le:r1r2}
If $r_1 \geq r_2 \geq r_3 \geq r_4 \geq 0$, and \eqref{eq:AB} holds,
then
\[
\frac{1}{r_i} + \frac{1}{r_j} \geq \frac{1}{r_k} + \frac{1}{r_l} - \frac{2}{r_k+r_l},
\]
where $i,j,k,l$ is any permutation of $\{1,2,3,4\}$.
\end{lemma}
\textit{Proof.}
Rewriting \eqref{eq:AB} as
\[
\left(\frac{1}{r_1} + \frac{1}{r_2} + \frac{1}{r_3+r_4} \right) \geq \left( \frac{1}{r_3} + \frac{1}{r_4} - \frac{1}{r_3+r_4} \right).
\]
Because $r_1 \geq r_2 \geq r_3 \geq r_4 \geq 0$, with some algebra one can verify that any other permutation will make the left-hand side of the preceding larger and the right-hand side of the preceding smaller, so the inequality still holds.
\hfill $\Box$

\begin{lemma}
\label{le:r4}
If $r_1 \geq r_2 \geq r_3 \geq r_4 \geq 0$, and \eqref{eq:AC} holds,
then
\[
\frac{1}{r_l} \leq \frac{1}{r_i+r_l} + \frac{1}{r_j+r_l} + \frac{1}{r_k+r_l},
\]
where $i,j,k,l$ is any permutation of $\{1,2,3,4\}$.
\end{lemma}
\textit{Proof.}
Multiplying by $r_4 (r_1+r_4) (r_2+r_4) (r_3+r_4)$ on both sides of \eqref{eq:AC} and canceling out common terms, we obtain
\[
r_1 r_2 r_3 \leq ((r_1 + r_2 + r_3 +r_4) + r_4) \, r_4^2.
\]
Because $r_4$ is the smallest, it is clear that any other permutation will make the left-hand side of the preceding smaller and the right-hand side of the preceding larger, so the inequality still holds.
\hfill $\Box$

\medskip

We are now ready for the main result in this subsection.

\begin{theorem}
\label{th:A}
Strategy A is optimal for the Searcher and the value of the game is $V_A$ if and only if
\begin{equation}
\frac{1}{r_1} + \frac{1}{r_2} \geq \frac{1}{r_3} + \frac{1}{r_4} - \frac{2}{r_3+r_4},
\label{eq:A>=B}
\end{equation}
and
\begin{equation}
\frac{1}{r_4} \leq \frac{1}{r_1+r_4} + \frac{1}{r_2+r_4} + \frac{1}{r_3+r_4},
\label{eq:A>=C}
\end{equation}
\end{theorem}
\textit{Proof.}
If strategy A is optimal for the Searcher, then $V_A \geq V_B$, which is equivalent to \eqref{eq:A>=B}, and $V_A \geq V_C$, which is equivalent to \eqref{eq:A>=C}.
Therefore, \eqref{eq:A>=B} and \eqref{eq:A>=C} are necessary conditions.

We next prove \eqref{eq:A>=B} and \eqref{eq:A>=C} are sufficient conditions.
Since the Searcher can use strategy A to guarantee an expected payoff at least $V_A$, it remains to show that the Hider has a mixed strategy to guarantee an expected payoff no more than $V_A$.
Let $q_{ij}$ denote the probability that the Hider hides the 2 booby traps in boxes $i$ and $j$, with probability $q_{ij} \geq 0$ and $\sum_{1 \leq i < j \leq 4} q_{ij} = 1$.
In particular, we will show that the Hider has a feasible mixed strategy to achieve an expected payoff exactly $V_A$ if the Searcher opens any one box, and guarantees an expected payoff no more than $V_A$ if the Searcher opens any two boxes.
In other words, we claim that there exists a feasible solution to
\begin{align*}
\sum_{1 \leq i < j \leq 4} q_{ij} &= 1, & \quad q_{ij} &\geq 0, \qquad \text{for } 1 \leq i < j \leq 4 \\
(q_{23} + q_{24} + q_{34} ) r_1 &= V_A, & \quad q_{12} (r_3 + r_4) &\leq V_A, \\
(q_{13} + q_{14} + q_{34} ) r_2 &= V_A, & \quad q_{13} (r_2 + r_4) &\leq V_A, \\
(q_{12} + q_{14} + q_{24} ) r_3 &= V_A, & \quad q_{14} (r_2 + r_3) &\leq V_A, \\
(q_{12} + q_{13} + q_{23} ) r_4 &= V_A, & \quad q_{23} (r_1 + r_4) &\leq V_A, \\
&& \quad q_{24} (r_1 + r_3) &\leq V_A, \\
&& \quad q_{34} (r_1 + r_2) &\leq V_A.
\end{align*}

To proceed, write $x = q_{34} / V_A$ and $y = q_{24} / V_A$, and use the first 5 equality constraints (in the left column) to solve $q_{ij}/V_A$ in terms of $x$ and $y$ for $1 \leq i < j \leq 4$.
Use $q_{ij} \geq 0$ to obtain lower bounds for $q_{ij}/V_A$, for $1 \leq i < j \leq 4$, and the last 6 inequality constraints (in the right column) to obtain their upper bounds.
The results are summarized below.
\begin{align*}
\frac{1}{(r_3 + r_4)} &\geq \frac{q_{12}}{V_A} = x + \frac{1}{2} \left( - \frac{1}{r_1} - \frac{1}{r_2} + \frac{1}{r_3} + \frac{1}{r_4} \right) \geq 0, \\
\frac{1}{(r_2 + r_4)} &\geq \frac{q_{13}}{V_A} = y + \frac{1}{2} \left( - \frac{1}{r_1} + \frac{1}{r_2} - \frac{1}{r_3} + \frac{1}{r_4} \right) \geq 0, \\
\frac{1}{(r_2 + r_3)} &\geq \frac{q_{14}}{V_A}  = \frac{1}{2} \left( \frac{1}{r_1} + \frac{1}{r_2} + \frac{1}{r_3} -\frac{1}{r_4} \right) - x - y \geq 0, \\
\frac{1}{(r_1 + r_4)} &\geq \frac{q_{23}}{V_A} = \frac{1}{r_1} - x - y \geq 0, \\
\frac{1}{(r_1 + r_3)} &\geq \frac{q_{24}}{V_A} = y \geq 0, \\
\frac{1}{(r_1 + r_2)} &\geq \frac{q_{34}}{V_A} = x \geq 0.
\end{align*}
Rewrite the preceding in terms of $x$, $y$, and $x+y$, to get the following.
\begin{align}
0 \leq x &\leq  \frac{1}{r_1+r_2},  \label{eq:x1} \\
- \frac{1}{2} \left( -\frac{1}{r_1} - \frac{1}{r_2} + \frac{1}{r_3} + \frac{1}{r_4} \right) \leq x &\leq \frac{1}{r_3+r_4} - \frac{1}{2} \left( -\frac{1}{r_1} - \frac{1}{r_2} + \frac{1}{r_3} + \frac{1}{r_4} \right), \label{eq:x2} \\
0 \leq y &\leq  \frac{1}{r_1+r_3}, \label{eq:y1} \\
- \frac{1}{2} \left( -\frac{1}{r_1} + \frac{1}{r_2} - \frac{1}{r_3} + \frac{1}{r_4} \right) \leq y &\leq \frac{1}{r_2+r_4} -  \frac{1}{2} \left( -\frac{1}{r_1} + \frac{1}{r_2} - \frac{1}{r_3} + \frac{1}{r_4} \right) \label{eq:y2}, \\
\frac{1}{r_1} - \frac{1}{r_1+r_4} \leq x + y &\leq \frac{1}{r_1}, \label{eq:xy1} \\
\frac{1}{2} \left( \frac{1}{r_1} + \frac{1}{r_2} + \frac{1}{r_3} -\frac{1}{r_4} \right) - \frac{1}{r_2+r_3} \leq x + y &\leq  \frac{1}{2} \left( \frac{1}{r_1} + \frac{1}{r_2} + \frac{1}{r_3} - \frac{1}{r_4} \right). \label{eq:xy2}
\end{align}
It then remains to show that there exists a feasible solution to $x$ and $y$ that satisfies these six linear constraints.

First, we claim there exists a feasible solution to $x$ that satisfies the two constraints \eqref{eq:x1} and \eqref{eq:x2}.
The larger lower bound for $x$ is clearly 0, since $r_1 \geq r_2 \geq r_3 \geq r_4$.
While it is not clear which of the two upper bounds for $x$ is smaller, one can verify that both are nonnegative, due to \eqref{eq:A>=B}.
With a similar argument, there exists a feasible solution to $y$ that satisfies \eqref{eq:y1} and \eqref{eq:y2}, due to \eqref{eq:A>=B} and Lemma~\ref{le:r1r2}.
 

Second, there exists a feasible solution to $x+y$ that satisfies \eqref{eq:xy1} and \eqref{eq:xy2}, because each of the two upper bounds is greater than or equal to each of the two lower bounds, due to  \eqref{eq:A>=B} and Lemma~\ref{le:r1r2}.

To complete the proof, we need to show that the sum between the upper bound (lower bound, respectively) of $x$ implied by \eqref{eq:x1} and \eqref{eq:x2} and the upper bound (lower bound, respectively) of $y$ implied by \eqref{eq:y1} and \eqref{eq:y2} is greater than or equal to the lower bound (upper bound, respectively) of $x+y$ implied by \eqref{eq:xy1} and \eqref{eq:xy2}.

From \eqref{eq:x1}, \eqref{eq:x2}, \eqref{eq:y1}, and \eqref{eq:y2}, the lower bound is 0 for $x$ and $y$, so we need to check the right-hand sides of \eqref{eq:xy1} and \eqref{eq:xy2} are both nonnegative.
The part concerning \eqref{eq:xy1} is trivial, and the part concerning \eqref{eq:xy2} follows because
\[
\frac{1}{r_1} + \frac{1}{r_2} + \frac{1}{r_3} \geq  \frac{1}{r_1+r_4} + \frac{1}{r_2+r_4} + \frac{1}{r_3+r_4} \geq  \frac{1}{r_4},
\]
where the second inequality follows from \eqref{eq:A>=C}.

Finally, we need to show that the sum between the upper bound of $x$ implied by \eqref{eq:x1} and \eqref{eq:x2} and the upper bound of $y$ implied by \eqref{eq:y1} and \eqref{eq:y2} is greater than or equal to the lower bound of $x+y$ implied by \eqref{eq:xy1} and \eqref{eq:xy2}.
We do so by showing that the sum of either upper bound of $x$ in \eqref{eq:x1} or \eqref{eq:x2}, and either upper bound of $y$ in \eqref{eq:y1} or \eqref{eq:y2}, is greater than or equal to either lower bound of $x+y$ in \eqref{eq:xy1} or \eqref{eq:xy2}.
There are thus 8 inequalities to verify.
For example, from \eqref{eq:x1}, \eqref{eq:y1}, \eqref{eq:xy1}, we need to show that
\[
\frac{1}{r_1+r_2} + \frac{1}{r_1+r_3} \geq \frac{1}{r_1} - \frac{1}{r_1+r_4},
\]
which follows from \eqref{eq:A>=C} and Lemma~\ref{le:r4}.
Using \eqref{eq:A>=C} and Lemma~\ref{le:r4}, we can also verify the corresponding inequality involving \eqref{eq:x1}, \eqref{eq:y2}, \eqref{eq:xy2}, and that involving \eqref{eq:x2}, \eqref{eq:y2}, \eqref{eq:xy2}, and that involving \eqref{eq:x2}, \eqref{eq:y1}, \eqref{eq:xy1}.

We next verify the corresponding inequality involving \eqref{eq:x1}, \eqref{eq:y1}, \eqref{eq:xy2}, which requires
\[
\frac{1}{r_1+r_2} + \frac{1}{r_1+r_3} \geq \frac{1}{2} \left( \frac{1}{r_1} + \frac{1}{r_2} + \frac{1}{r_3} -\frac{1}{r_4} \right) - \frac{1}{r_2+r_3},
\]
which is equivalent to
\begin{equation}
\frac{2}{r_1+r_2} + \frac{2}{r_1+r_3} \geq  \left( \frac{1}{r_1} + \frac{1}{r_2} + \frac{1}{r_3} -\frac{1}{r_4} \right) - \frac{2}{r_2+r_3}.
\label{eq:r2r3}
\end{equation}
Starting from the right-hand side to get
\begin{align*}
\left( \frac{1}{r_2} + \frac{1}{r_3} - \frac{2}{r_2+r_3} \right) + \frac{1}{r_1} - \frac{1}{r_4} &\leq \left( \frac{1}{r_1} + \frac{1}{r_4} \right) + \frac{1}{r_1} - \frac{1}{r_4} \\
&= \frac{2}{r_1+r_1} + \frac{2}{r_1+r_1} \\
&\leq \frac{2}{r_1+r_2} + \frac{2}{r_1+r_3},
\end{align*}
where the first inequality follows from \eqref{eq:A>=B} and Lemma~\ref{le:r1r2}, and the last inequality follows from $r_1 \geq r_2$ and $r_1 \geq r_3$.

We can go through the same procedure to  verify the corresponding inequality involving \eqref{eq:x1}, \eqref{eq:y2}, \eqref{eq:xy1}, and that involving \eqref{eq:x2}, \eqref{eq:y1}, \eqref{eq:xy1}, and that involving \eqref{eq:x2}, \eqref{eq:y2}, \eqref{eq:xy2}.
Each of these three inequalities has the same form as in \eqref{eq:r2r3}, with the bracket on the right-hand side having three positive terms and 1 negative term.
The key to establish the inequality is to apply Lemma~\ref{le:r1r2} to the two positive terms with the largest indices among the three positive terms; for example, to prove \eqref{eq:r2r3} we pick $r_2$ and $r_3$ to apply Lemma~\ref{le:r1r2}.

Because there exists feasible solution to $x$ and $y$ that satisfies constraints \eqref{eq:x1}--\eqref{eq:xy2}, we have shown that the Hider has a mixed strategy that guarantees an expected payoff no more than $V_A$, which completes the proof.
\hfill $\Box$

\subsection{General bounds}
Here we give some general bounds on the value of the game, starting with an upper bound and a lower bound that are close to each other when $n$ is large and the rewards are small.

\begin{proposition}
	(Upper bound)
	The value $V$ of the game on the complete hypergraph satisfies
	\begin{align}
	V \le  \frac{R_0}{k+1}\left(1- \frac{1}{k+1} \right)^k, \label{eq:Vub}
	\end{align}
	where $R_0 = \sum_{i=1}^n r_i$.
\end{proposition}
\textit{Proof.}
First assume all the rewards are integers. If we define a new game by replacing a box of reward $r$ with two boxes of reward $r_1$ and $r_2$ with $r_1+r_2 = r$, then the value of the game can only increase, because any Searcher strategy in the original game can also be used in the new game.
	With a similar argument, we can replace each box $i$ with $r_i$ boxes each containing a reward of $1$, resulting in a new game with equal rewards of 1, whose value is at least as great as the original game.
	The value of the new game is equal to $U(R_0, k)$, as defined in Theorem~\ref{thm:eq}.
	Observe that by further replacing each box with $t$ new boxes each containing a reward of $1/t$, we obtain a game whose value $U(t R_0, k)/t$ is no smaller than that of the original game.
	Therefore, the value of the original game is bounded above by
\[
\lim_{t \rightarrow \infty} \frac{U(t R_0, k)}{t} = R_0 \lim_{t \rightarrow \infty} \frac{U(t R_0, k)}{t R_0} = \frac{R_0}{k+1}\left(1- \frac{1}{k+1} \right)^k,
\]
where the last equality follows from \eqref{eq1}.
	
	If the rewards are all rational numbers, then we can obtain an equivalent game with integer rewards by multiplying them all by a common denominator $d$. 
	All the payoffs in the resulting game will be larger by a factor of $d$, and therefore so will the value of the game and the parameter $R_0$.
	As a consequence, the left- and right-hand sides of~(\ref{eq:Vub}) will both be larger by a factor of $d$, so the inequality still holds.
	If the rewards are real numbers, then they can be approximated arbitrarily closely to rational numbers, so that the left- and right-hand sides of~(\ref{eq:Vub}) are also approximated arbitrarily closely, and the bound still holds.
\hfill $\Box$

\bigskip

\begin{proposition} \label{prop:indep}
	(Lower bound)
	The value $V$ of the game on the complete hypergraph satisfies
	\begin{align}
	V \ge  \frac{R_0}{k+1}\left(1- \frac{1}{k+1} \right)^k \left( 1 - \frac{r([k])}{R_0} \right), \label{eq:Vlb}
	\end{align}
	where $R_0 = \sum_{i=1}^n r_i$ and $r([k]) = \sum_{i=1}^k r_i$.  
\end{proposition}
\textit{Proof.}
Consider a Searcher strategy with which each box is independently opened with probability $1/(k+1)$.
For a given Hider strategy $H \in [n]^{(k)}$, the probability that none of the boxes in $H$ is opened is $(1-1/(k+1))^k$.
If the Searcher does not open any box in $H$, her expected payoff is $r(\bar{H})/(k+1)$; if she opens any boxes in $H$, her payoff is zero. 
Therefore, with such strategy the Searcher's expected payoff is
\[
\left(1- \frac{1}{k+1} \right)^k \left( \frac{r(\bar H)}{k+1} \right).
\]
The preceding in minimized when $r(H)$ is maximized; that is, for $H = [k]$. In this case, the expected payoff is the right-hand side of~(\ref{eq:Vlb}).
\hfill $\Box$

\bigskip

It is worth pointing out that, among all the Searcher strategies that open each box independently at random with some given probability $p$, the one that guarantees the greatest expected payoff is given by $p=1/(k+1)$, namely the strategy of Proposition~\ref{prop:indep}.
This claim can be verified via elementary calculus.
The bounds in~(\ref{eq:Vub}) and~(\ref{eq:Vlb}) are close when $r([k])/R_0$ is close to zero. In particular, the bounds are asymptotically equal for constant $k$, as $n \rightarrow \infty$, if all the rewards are all $o(n)$. In this case, the Searcher strategy that opens each box independently with probability $1/(k+1)$ is asymptotically optimal. 

Note that all the optimal Searcher strategies presented in this paper share the same form: the Searcher chooses each hyperedge $S$ with probability 0, or with probability proportional to $1/r(S)$. This observation gives rise to a set of lower bounds on the value, generalizing the Searcher strategy from Lemma~\ref{lem:null}. 

\begin{proposition} \label{prop:partition}
Consider the search game played on an arbitrary hypergraph, and let $\mathcal{S} = \{S_1, \ldots, S_t \}$ be a set of hyperedges. Consider the Searcher strategy $p$ that chooses $S_j$ with probability $p(S_j) = \lambda/r(S_j)$, where
	\[
	\lambda \equiv  \lambda(\mathcal S) \equiv \frac{1}{\sum_{j=1}^t 1/r(S_j)}.
	\]
	This strategy guarantees an expected payoff of at least $m \lambda$, where
	\[
	m \equiv m(\mathcal S) \equiv \min_{H \in [n]^{(k)}} |\{S_j \in \mathcal{S}: S_j \cap H = \emptyset \} |
	\]
	is the minimal---over all possible Hider strategies---number of hyperedges in $\mathcal S$ that contain no booby traps.
\end{proposition}
\textit{Proof.}
For a given Hider strategy $H$, let $ \mathcal{A} = \{S_j \in \mathcal{S}: S_j \cap H = \emptyset \} |$ be the set of hyperedges in $\mathcal{S}$ that contain no booby traps. By definition of $m$, we have $|\mathcal{A} | \ge m$. Hence, the expected payoff against $H$ is
\[
R(p,H) = \sum_{S \in \mathcal A} p(S) r(S) = \sum_{S \in \mathcal A} \lambda \ge m \lambda,
\]
which completes the proof.
\hfill $\Box$

If $\mathcal{S}$ is a partition of $[n]$, then the minimal number of hyepredges that contain no booby traps is $m(\mathcal{S})=t-k$, and Proposition~\ref{prop:partition} implies that the value is at least $(t-k) \lambda(\mathcal{S})$.

Based on the solutions to special cases presented in this paper, we make a conjecture on the Searcher's optimal strategy.

\begin{conjecture}
Consider the booby trap search game played on a hypergraph.
There exists an optimal Searcher strategy with which each hyperedge will not be chosen at all, or will be chosen with probability inversely proportional to the sum of the rewards on that hyperedge.
In other words, the Searcher can achieve optimality by choosing the best subset of hyperedges and using the mixed strategy described in Proposition~\ref{prop:partition}.
\end{conjecture}

\section{Conclusion}
\label{sec:conclusion}
This paper presents a new search game on a hypergraph between a Searcher and a Hider.
The Searcher wants to collect maximum reward but has to avoid booby traps planted by the Hider.
We present the solutions to a few special cases, based on which we make a conjecture about the form of the solution in general.

Two of the special cases presented in this paper involve the Searcher opening just one box, or opening any number of boxes.
A relevant and practical situation may restrict the Searcher to opening a certain fixed number of boxes.
If the booby trap only partially injures the Searcher but does not incapacitate her, then we can consider a model extension that allows the Searcher to keep going until she encounters a certain number of booby traps.

\section*{Acknowledgements} This material is based upon work supported by the National Science Foundation under Grant No. IIS-1909446.

\bibliographystyle{apalike}
\bibliography{../bib/references}

\end{document}